\documentstyle[12pt]{article}
\newfont{\mib}{cmmib10 scaled 1200}

\renewcommand{\phi}{\varphi}
\newtheorem{Lemma}{Lemma}[section]
\newtheorem{Theorem}{Theorem}[section]
\newtheorem{Proposition}{Proposition}[section]

\newtheorem{Corollary}[Proposition]{Corollary}
\newtheorem{Definition}{Definition}[section]

\newtheorem{Remark}{Remark}[section]
\newtheorem{Conjecture}{Conjecture}[section]
\begin{document}
\title{Minimal discrepancies of toric singularities}
\author{Alexandr Borisov\\
Department of Mathematics \\
Pennsylvania State University\\
e-mail: borisov@math.psu.edu}
\date{July 7, 1994}
\maketitle

\begin{abstract}

The main purpose of this paper is to prove that minimal discrepancies of
$n$-dimensional toric singularities can accumulate only from above and only
to minimal discrepancies of toric singularities of dimension less than $n$. I
also prove that some lower-dimensional minimal discrepancies do appear as such
limit.
\end{abstract}

\section{Introduction}
First of all let me stress that minimal discrepancy of a variety $X$ I am
talking about is defined only for singular varieties and is the minimum of
discrepancies of exceptional divisors in all resolutions of singularities of
$X$. I don't allow to blow up anything in a smooth part of $X$ otherwise the
minimal discrepancy would always be no more than $1$ and this means loss of
information about $X.$

Now let me recall the following conjecture proposed by V. Shokurov.(\cite{Sh})

\begin{Conjecture}
For every natural $n$ minimal discrepancies of $n$-dimensional log-terminal
singularities can accumulate only from above.
\end{Conjecture}

In particular this conjecture implies that for every $n$ there exists a
positive constant $\epsilon (n),$ such that if all discrepancies of
$n$-dimensional variety $X$ are greater than $-\epsilon (n)$ then they are in
fact nonnegative, that is $X$ has at most canonical singularities.

In this paper the above conjecture will be proven for a particular case of
toric singularities. By the simple observation (see Lemma $2.1$) the problem
can be reduced to the case of cyclic quotient singularities for which some
results are obtained. They are in some sense best possible for nonterminal
singularities and also very informative for terminal case. What to do in case
of more general singularities is discussed in section $3.$ The main result of
the paper is the following. (Corollary 2.1)

{\bf Main Result. For every natural $n$  minimal discrepancies of
$n$-dimensional toric singularities can accumulate only from above and only
to minimal discrepancies of toric singularities of dimension less than $n.$}

\begin{Remark}
It can happen that infinitely many different toric singularities have the same
minimal discrepancy. I do not consider this as an accumulation of minimal
discrepancies. In this case the minimal discrepancy in question does not
necessarily come from lower dimension. It may be of the form (1+(minimal
discrepancy of lower-dimensional toric singularity.) However the only example
of that kind I know is rather trivial case of minimal discrepancy $n-1$ for
$2n-$dimensional singularities.
\end{Remark}

There is also the result in the opposite direction (see Theorem $2.2$) which
implies in particular that every minimal discrepancy of a toric singularity of
dimension $k$ is a limit of minimal discrepancies of $n-$dimensional toric
singularities for $n$ big enough. And it also implies that if this discrepancy
is non-positive then the only restriction on $n$ is that $n>k.$

I am glad to thank here V. Shokurov for his constant interest in this study and
helpful remarks about it.


\section{Proofs}
First of all, let me recall some basic facts about toric varieties. (See, for
example, \cite {BB}, \cite {D}, \cite {R}.) Every $n$-dimensional affine toric
variety $X$ is just a $Spec(R),$ where R is a ring generated by monomials
$x_1^{\alpha _1} x_2^{\alpha _2}...x_n^{\alpha _n},$ where $\{\alpha _1,\alpha
_2,...,\alpha _n\}$ is an integral point of some finitely generated convex cone
of full dimension $C(X)$ in $R^n.$ For several reasons it is more useful to
consider the dual cone $C^*(X)$ in the dual space $V=R^n.$ It does not
necessarily have dimension $n,$ but we will always assume that it will.
Otherwise $X$ would be isomorphic (not canonically) to a product of another
toric variety and an algebraic torus of positive dimension.

Now various conditions on singularities of $X$ have simple combinatorial
formulation in terms of $C^*(X).$ Namely, let us consider one-dimensional
extremal rays $l_1,l_2,...,l_k$ that generate this cone. They are rational
therefore we can pick on every $l_i$ an integral point $P_i$ which is the
closest one to zero. Then all $P_i$ lie in one hyper-plane if and only if $X$
is Q-Gorensteinian. $X$ is Q-factorial if and only if $k=n,$ that is $C^*$ (or
$C$) is simplicial. $X$ is regular if and only if $C^*(X)$ is regular which
means that it is simplicial and $P_i$ form a basis for the lattice. Moreover,
the Gorenstein index and minimal discrepancy of $X$ also have simple
description. Namely, let us consider the linear function $F$ on $V,$ such that
$F(P_i)=1.$ (This is possible exactly when $X$ is Q-Gorensteinian.) Then the
least common denominator of values of $F$ on non-zero points of $C^*$ is the
index. The minimal log-discrepancy which is by definition (1+(minimal
discrepanc!
y)) is the minimum of the above values among points in the interior of
non-regular sub-cones of $C^*(X).$ (If all sub-cones including $C^*(X)$ itself
are regular then $X$ is regular and its minimal discrepancy is undefined.) We
will pass freely from discrepancies to log-discrepancies mostly using the
latter in proofs and the former in statements.

There is one type of toric singularities which is particularly interesting for
our purposes. Namely, quotients of an affine plane $A^n$ by cyclic groups. This
corresponds to the case when $C^*$ is simplicial and the lattice $N$ of all
integral points in $V$ is generated by $P_i$ and only one extra element $x.$ We
can always assume that $x$ lies in the interior of $C^*$ otherwise the
singularity splits into the lower-dimensional singularity and torus. I want to
mention that the cyclic group in question is just a factor ($N/<P_i>$) and its
action can be reconstructed from the coordinates of $x$ in the basis $\{P_i\}.$
The following lemma reduces everything to this special case.
\begin{Lemma} The set of minimal discrepancies of toric singularities of
dimension $n$ coincides with that of cyclic quotients of dimension no greater
than $n.$
\end{Lemma}
{\bf Proof} Let $X$, $C^*\subset V$, $P_i$, $F$ be as above. Let $\epsilon$ be
the minimal log-discrepancy of $X.$ By the above combinatorial description
there exists an integral point $x\in C^*,$ such that $F(x)=\epsilon .$ Consider
the ray generated by $x.$ It intersects the polygon $P_1P_2...P_k$ in some
point $P.$ Standard combinatorial arguments show that there exists a simplex
$P_{i_1}P_{i_2}...P_{i_r}\subset P_1P_2...P_k,$ such that $P$ lies in its
interior. Of course, this simplex has dimension no greater than $n-1$ and its
interior means interior with respect to its own geometry.

Now let us stick to the subspace $W$ of $V$ generated by $P_{i_j}.$ Evidently,
$C^*\bigcap W$ will be a convex cone corresponding to some new toric variety
$X^\prime$ with the same minimal log-discrepancy $\epsilon .$ This $X^\prime$
is already Q-factorial but it is not a cyclic quotient yet. To produce out of
it a cyclic quotient let me do the following. Consider lattice $N^\prime
\subset N$ generated by $x,$ change coordinates in such a way that $N^\prime$
become a lattice of integral points and forget about $N.$ What we have now is a
cyclic quotient $X^{\prime \prime},$ which again has the same minimal
log-discrepancy $\epsilon,$ so the lemma is proven.

{\bf Remark} As one can easily see from the proof of the above lemma we can
assume that factor-group ($N/<P_i>$) is not only cyclic but also generated by
the element $x,$ which has "minimal log-discrepancy" (that is $F(x)=\epsilon
.$)
In the rest of the paper this element will be often called generating element.
The fact that it is not uniquely determined for a given singularity does not
cause any difficulties.

{}From now on we stick to this particular case of quotient singularities and
whenever we have a toric variety it is a cyclic quotient singularity. The above
lemma allows us to do it. Now let me notice that the results we are going to
prove are of two types. Most of them are negative in a sense that they restrict
where and how minimal discrepancies can accumulate. And there are some positive
results based on procedures that allow us to construct cyclic quotients with
prescribed minimal discrepancies starting with the given one. We begin with
negative results which all deal with the following situation.

 Suppose we have a sequence of cyclic quotient singularities $\{X^\nu\}, \nu
=1,2,...,$ such that their log-discrepancies $\epsilon ^\nu$ are getting closer
and closer to some real number $\epsilon. $ Consider the standard simplex
$\Delta$ in $R^n$ defined by the inequalities $\alpha_i \ge 0$, $\sum
{\alpha_i} \le 1$ inside a standard hypercube $H$, defined by the inequalities
$0\le \alpha_i \le 1.$ By identifying simplexes $P_1^\nu P_2^\nu ...P_k^\nu $
with this standard one we obtain a sequence of points $\alpha^\nu \in H$ that
correspond to $x^\nu .$ By the compactness of $H$ there exists a subsequence
with a limit point $\alpha.$ We replace our sequence by this subsequence. By
the above remark $\epsilon ^\nu = \sum {\alpha_i^\nu}.$ Therefore $\epsilon =
\sum {\alpha_i}.$
Now all our negative results can be formulated as the sequence of statements
which will be proven in a row. Let me state this as a theorem.
\begin{Theorem} In the above notations the following is true.
Content-Length: 17452

$1)$ If $\epsilon ^\nu$ are not the same for big $\nu$ then $\alpha$ is on the
boundary of $H.$

$2)$ If $\epsilon ^\nu$ are the same for big $\nu$ then one can choose $\alpha$
on the boundary of $H$ which has the same $\epsilon$ and is also a limit of
some sequence of the same type.

$3)$ $\epsilon$ is rational.

$4)$ If $\epsilon^\nu$ are not the same for big $\nu$ they accumulate to
$\epsilon$ only from above.

$5)$ We can choose $\alpha$ as in $2)$ on some face of $H$ to be a generating
point of a cyclic quotient singularity if considered on this face. As a
corollary, $\epsilon$ is a minimal log-discrepancy for some lower-dimensional
toric singularity plus some nonnegative integer.

$6)$ Every face of $H$ is characterized by restricting some coordinates to be
$0$ and some coordinates to be $1.$ Under this remark statement 5 can be
strengthen by the restriction that for the face of $\alpha$ the number of $1-$s
is not greater than the number of $0-$s. Moreover if $\epsilon ^\nu$ are not
the same for big $\nu$ then the  number of $1-$s is strictly less than the
number of $0-$s
\end{Theorem}
In order to prove this theorem let me introduce the notion of multiple of the
point in $H.$ It will be used a lot in the rest of the paper so it deserves to
be stated formally.

\begin{Definition}
For every point $\alpha =(\alpha_i)\in H$ and integer $m$ let $m-$th multiple
of $\alpha$ be the point $\alpha^{(m)}$ whose $i-th$ coordinate is $1$ if
$\alpha_i=1$ and $\{m\alpha_i\}$ otherwise. Note that for positive $m$ this
construction is continuous at the neighborhood of the boundary of $H.$
\end{Definition}

Now we begin the proof.

{\bf Statement 1).} Suppose $\alpha$ is in the interior of $H.$ Consider
$\alpha^{(m)}$ for all integer $m.$ Then the compactness of $H$ tells us that
there are two numbers $m_1 < m_2,$ such that $\alpha^{(m_i)}$ are very close,
for example closer than $\frac1{100} \times$ (distance from $\alpha$ to the
boundary of H). They may also coincide, we don't care. Then
$\alpha^{(m_1-m_2+1)}$ and $\alpha^{(m_2-m_1+1)}$ are evidently very close to
$\alpha^{(1)} = \alpha.$ We have several cases.

First of all suppose that sum of coordinates of one of the above two points is
less than $\epsilon$ (that means that sums of coordinates of
$\alpha^{(m_1-m_2+1)}$ and $\alpha^{(m_2-m_1+1)}$ are different.) Let it be
$\alpha^{(m_1-m_2+1)}.$ Then for $\nu$ big enough $\alpha^\nu$ is close enough
to $\alpha$ and $\alpha^{\nu,(m_1-m_2+1)}$ is close enough to
$\alpha^{(m_1-m_2+1)}$ and therefore the sum of coordinates of
$\alpha^{\nu,(m_1-m_2+1)}$ is less than $\epsilon^\nu,$ which is impossible.

Now suppose that sums of coordinates are the same. Then if there is a
subsequence of $\alpha^\nu$ for which $\epsilon^\nu$ accumulate to $\epsilon$
from above consider $(m_1-m_2+1)-$th multiples. Then for $\nu$ big enough from
this subsequence sum of the coordinates of $\alpha^{\nu,(m_1-m_2+1)}$ is less
then sum of the coordinates of $\alpha^{(m_1-m_2+1)},$ because $(m_1-m_2+1)
<0.$ Therefore it is less then  $\epsilon^\nu,$ which is impossible. Similar
arguments work for the case when $\epsilon^\nu$ accumulate to $\epsilon$ from
below. We should just consider $(m_2-m_1+1)-$th multiples instead of
$(m_1-m_2+1)-$th ones and notice that $(m_2-m_1+1) \ge 2.$

\begin{Remark}
We did not prove that point inside $H$ cannot be a limit of generating points
of cyclic quotients with the {\bf same} discrepancy. And this indeed can
happen. The easiest example is given by two-dimensional canonical toric
singularities.
\end{Remark}

{\bf Statement 2).} Suppose we have a sequence of points $\alpha^\nu$ with the
same sum of coordinates $\epsilon.$ Consider those multiples of all these
points that have the same sum of coordinates $\epsilon.$ We will see very soon
that there are plenty of them. We have two cases.

First of all, suppose $\alpha$ has finite order in $H$ that is $\alpha^{(k)} =
0$ for some $k.$ Then for $\alpha^\nu$ that are close enough to $\alpha$
$(mk+1)-$multiples have sum of the coordinates $\epsilon$ Moreover when we make
$m$ run from zero to some number depending on $\nu$ they run following some
straight line with small intervals until they hit the boundary of $H.$ The
length of these intervals goes to zero when $\alpha^\nu$ go to $\alpha.$
Therefore we have infinitely many points in every neighborhood of the boundary
of $H,$ intersected with a hyper-plane $\sum{x_i} =\epsilon.$ Therefore there
exists a point on this boundary which is a limit of some sequence of these
points. To complete the argument it is enough to mention that each one of these
points is also a generating point for some quotient singularity with the same
discrepancy $\epsilon.$

Now suppose $\alpha$ has infinite order in $H.$ Nevertheless $\epsilon$ is
rational because $\epsilon =\epsilon^\nu .$ So we have infinitely many
multiples of $\alpha$ with the same sum of coordinates $\epsilon.$ Then
arguments similar to that of the above case show that whenever two of such
multiples are close to each other there is some other multiple that is close to
the boundary of $H.$ Again as before, there is a point on the boundary which is
a limit of a sequence of these multiples. Now we can just notice that every
multiple of $\alpha$ is a limit of multiples of $\alpha^\nu$ and the rest is
the same as above.

{\bf Statement 3).} Suppose $\epsilon$ is irrational. By previous statements we
can assume that $\alpha$ is on some face of $H.$ Then all its multiples are by
the definition on the same face. Now we want to prove that for some $m>0$
$\alpha ^{(m)}$ is close enough to $\alpha$ and sum of the coordinates of
$\alpha ^{(m)}$ is less than $\epsilon .$ This is not completely trivial,
because we require $m$ to be positive. Here is the proof. First of all, we can
stick to the face of $H$ $\alpha$ belongs to. Then we notice that all sums of
coordinates of $\alpha ^{(m)}$ are different because $\epsilon$ is irrational.
By the compactness argument there exist some positive integers $m_1<m_2$ such
that $\alpha ^{(m_1)}$ and $\alpha ^{(m_2)}$ are close enough. If the sum of
coordinates of $\alpha ^{(m_1)}$ is greater than the sum of coordinates of
$\alpha ^{(m_2)}$ it is enough to choose $m$ to be equal to $1+m_2-m_1.$
Otherwise we need one more step. Denote $m_3=1+m_1-m_2.$ Then everyth!
ing would have been OK, but $m_3$ is not positive. But we can find $m_4 < m_5$
of form $l(m_3+1),$ such that $\alpha ^{(m_4)}$ and $\alpha ^{(m_5)}$ are so
close that $\alpha ^{(m_3+m_5-m_4)}$ is still close enough to $\alpha$ and the
sum of its coordinates is still less than the sum of coordinates of $\alpha.$

Now for $m$ as above and $\nu$ big enough $\alpha ^\nu$ is close enough to
$\alpha$ therefore $\alpha ^{\nu ,(m)}$ is close enough to $\alpha ^{(m)}.$
(Here we really need that $m$ is positive because $\alpha$ lies on the boundary
of $H.$) But this means that for $\nu$ big enough sum of the coordinates of
$\alpha ^{\nu ,(m)}$ is less than sum of the coordinates of $\alpha ^\nu,$
which is impossible.

{\bf Statement 4).} Now $\epsilon$ is rational. The arguments similar to the
above allow us to find an integer $m>1$ such that $\alpha ^{(m)}$ is close to
$\alpha $ and has the same sum of coordinates. Namely, the compactness argument
tells that there are $m_1,m_2$ such that $\alpha ^{\nu ,(m_1)}$ and $\alpha
^{\nu ,(m_2)}$ are arbitrary close. (They may even coincide, we don't care.)
Then $m=1+m_2-m_1$ will satisfy all requirements.

Now if discrepancies $\epsilon ^\nu$ accumulate to $\epsilon$ from below then
for sufficiently large $\nu$ the sum of coordinates of $\alpha ^{\nu ,(m)}$ is
less than the sum of coordinates of $\alpha ^\nu.$ (By definition $m$ is
greater than 1 and $(\alpha ^{\nu ,(m)}-\alpha ^{(m)})=m(\alpha ^\nu -\alpha
).$) This completes the proof of the statement.

{\bf Statement 5).} We have $\alpha$ on some face of $H.$ Let us consider this
face and multiples of $\alpha$ on it. If there are infinitely many of them that
are in fact different then there are infinitely many of them with the same sum
of the coordinates $\epsilon,$ because $\epsilon$ is rational. Then arguments
of the proof of statement $2$ allow us to replace $\alpha$ so that it lie on a
face of lower dimension. We can do this until we come to $\alpha$ that has
finite order in the appropriate face. Now on this face $\alpha$ is a generating
point for a quotient singularity, because if some multiple of it has smaller
sum of coordinates in the face it has smaller sum of all coordinates and usual
arguments show that it is impossible.

{\bf Statement 6).} Suppose $\alpha$ has order $N$ in its face. Suppose this
face is determined by $k$ equalities of type $x_i =0$ and $l$ equalities of
type $x_i =1$ Consider $(1-N)-$th multiples of $\alpha ^\nu .$ Then the
corresponding $\epsilon-$s go to $\epsilon +k-l$ when $\alpha ^\nu$ go to
$\alpha.$ Therefore $k \ge l.$ Moreover, if $\epsilon ^\nu$ accumulate to
$\epsilon$ from above then $\epsilon-$s for $(1-N)-th$ multiples accumulate
from below. So case $k=l$ is also impossible.

This completes the proof of the theorem. The following corollary is formally
also restrictive but in fact as you can see from its proof it is a positive
result (or, more precisely, simple observation.)

\begin{Corollary}
Under the notations of the above theorem if $\epsilon ^\nu$ are not the sam for
big $\nu$ then $\epsilon$ is not just sum of lower-dimensional log-discrepancy
and integer but is a lower-dimensional log-discrepancy itself. If $\epsilon
^\nu$ are the same for big $\nu$ then $\epsilon$ is either a lower-dimensional
log-discrepancy or (1+(minimal log-discrepancy of dimension $\le (n-2)$)).
\end{Corollary}

{\bf Proof} This is a straightforward consequence of statement $6$ and the
following fact.

{\bf Fact} For arbitrary $m-$dimensional cyclic quotient one can construct
$(m+2)-$dimensional cyclic quotient whose minimal log-discrepancy is greater
than given exactly by $1.$

Construction that proves the above fact is as follows. Suppose the generating
point $\alpha$ has order $N.$ Then we construct new $(m+2)-$dimensional
singularity defined by the generating point $(\alpha ,\frac1N, 1-\frac1N ) $
which means that first $m$ coordinates remain the same and last two are as
specified.

\begin{Remark}
The main idea of the proof of the above theorem (namely use of multiples,
compactness and some sort of continuity) is very similar to that of the
boundedness theorem for toric Fano varieties with bounded discrepancies. (\cite
{BB}) However there everything is written using less geometrical language that
maybe hide this idea in formulas. Probably it can also be written using the
language of this paper, I just don't know anybody who really tried. I can do it
for canonical singularities, but not in general case.
\end{Remark}

Now let me state the most general positive result I know about what
lower-dimensional discrepancies can indeed appear as a limit. But before doing
this I would like to notice that by the evident reason of symmetry every
minimal log-discrepancy of cyclic quotient of dimension $n$ is no greater than
$\frac n2 .$

\begin{Theorem} Suppose we have an $m-$dimensional cyclic quotient generated by
$\alpha$ with minimal discrepancy $\epsilon.$ Denote $r=-[-\epsilon]$ so that
$r$ is the smallest real number which is greater or equal than $\epsilon .$
Then for all nonnegative integers $l$, $(\epsilon +l)$ is a limit of
$n-$dimensional log-discrepancies for all $n \ge m+r+2l.$
\end{Theorem}

{\bf Proof} There are in fact several ways of doing this for nonzero $l.$ The
freedom we have is basically due to the Fact in the proof of the above
corollary. I will show you just one way.

First of all let me consider the standard $n-$dimensional hypercube $H$ and
divide the set of coordinates $\{x_1, x_2 ... x_n\}$ into three parts as
follows. First $m$ of them will correspond to the coordinates of our given
$m-$dimensional singularity and will be still called $x_i$. Those with indexes
from $m+1$ to $m+l$ will be called $y_i, i=1,...,l.$ And those with indexes
from $m+1+l$ to $n$ will be called $z_i, i=1,...,n-m-l.$ Now we place our
$m-$dimensional singularity on the face of $H$ defined by equalities
$y_i=1,z_i=0.$ Let us denote by $T$ the point $(\alpha;1,...,1;0,...,0)$ that
corresponds to the generating point $\alpha.$ Suppose its order is $q$ that is
its $q-$th multiple is a vertex of $H.$ Consider vertex
$P=(0,...,0;0,...,0;1,..,1).$ (Here ";" divides $x_i$, $y_i$ and $z_i.$) Now
the generating points of $n-$dimensional singularities we are looking for lie
on the segment $PT$ close to $T.$ More precisely, they are given by formula
$\frac1NP+(1-\frac1N)T)!
$ where $q|(N-1).$ In order to complete the proof it is enough to show that for
all such points every multiple is either trivial or has greater or equal sum of
coordinates than the point itself.

So, consider the point $A=A_N$ defined as above. It is a straightforward
observation that it has order $N.$ Now arguments similar to those from the
proof of the statement $2$ of Theorem $2.1$ show that for every positive
integer $k<N$ $k$-th multiple of $A$ lies in the set $S_k$ defined by the
following procedure.

{\bf Procedure.}
 Suppose $T_k$ is a $k-th$ multiple of $T.$ Draw a raw starting from $T_k$ and
parallel to the raw $[TP).$ When it hits the boundary of $H$ change the
corresponding $1-$s to $0-$s and $0-$s to $1-$s. Do it until the sum of lengths
of all segments drawn equals the length of $PT.$

 So it is enough to prove that no point from this set can have greater sum of
coordinates than that of $A.$ In order to do it let me make several simple
observations. First of all let me notice that in fact when we draw our segments
we $\bf never$ change $1$ to $0$, we only change $O-$s to $1-$s. The reason is
that in all $x_i$ we draw in the negative direction and we cannot hit the
boundary on $y$ or $z$ for $k<N.$  Another observation is that "locally", that
is when we don't hit boundary, the sum of coordinates does not decrease because
sum of coordinates of $P$ is greater or equal than that of $T$ by the condition
$n \ge m+r+2l.$ Combined together these two observations evidently take care of
$k$ which are not divisible by $q.$ For $k$ divisible by $q$ we only need to
notice that we begin our procedure from the point $(0,...,0;0,...,0;0,...,0)$
but the first nontrivial segment starts from the point
$(1,...,1;0,...,0;0,...,0)$ whose sum of coordinates is greater than that !
of $A.$

\section{Some open questions}
There are several natural questions concerning the obtained results.

{\bf Question 1.} Is it true that EVERY minimal log-discrepancy of
$n-$dimensional cyclic quotient is a limit of minimal log-discrepancies of
$(n+1)-$dimensional cyclic quotients? The theorem above together with the
classification of $3-$dimensional terminal toric singularities implies that
this is true for $n\le 3.$ It is natural to try, maybe with computer, the case
$n=4.$ As far as I know the classification of $4-$dimensional toric terminal
singularities is not yet completed but there are some conjectures and a lot of
work is already done. See \cite{MMM} for details.

Other questions naturally arose when I tried to extend these results to more
general singularities.

{\bf Question 2.} Is it true that the set of minimal discrepancies of quotient
singularities with respect to arbitrary groups coincides with that of cyclic
quotients of the same dimension?

{\bf Question 3.} Is there an example of log-terminal singularity whose minimal
discrepancy is not a minimal discrepancy for any cyclic quotient of the same
dimension?

{\bf Question 4.} Is it true for arbitrary log-terminal singularities which are
not terminal that every (or at least one) divisorial valuation that corresponds
to the minimal discrepancy is given by a divisor on the Q-factorial terminal
modification in sense of Clemens-Koll\'ar-Mori? (Of course it is not true for
all valuations with negative discrepancy, but the question is about minimal
discrepancy.) It is true for toric singularities and in dimension $2$ and I
have no counterexamples in general case.

While stating these questions it would be unfair not to express my opinion
about them. I suspect that the answer to Question $1$ is "Yes" for many
singularities but not for all of them. The answer to Question $2$ is probably
"Yes". Example to Question $3$ probably also exists maybe even $3-$dimensional.
And I have been unable so far to find any evidence pro or against to Question
$4.$

\end{document}